\documentclass{amsart} 
\usepackage{amssymb,latexsym,amscd, url}

\theoremstyle{plain}
\newtheorem{theorem}{Theorem}[section]
\newtheorem{corollary}[theorem]{Corollary}
\newtheorem{lemma}[theorem]{Lemma}
\newtheorem{proposition}[theorem]{Proposition}
\newtheorem{alg}[theorem]{Algorithm}
\theoremstyle{definition}

\newtheorem{remark}[theorem]{Remark}

\newcommand {\Q}{{\mathbb{Q}}}
\newcommand {\C}{{\mathbb{C}}}
\newcommand {\R}{{\mathbb{R}}}
\newcommand {\Z}{{\mathbb{Z}}}
\newcommand {\N}{{\mathbb{N}}}
\newcommand {\F}{{\mathbb{F}}}

\newcommand {\OO}{{\mathcal{O}}}

\newcommand{\POL}  {\mathop{\rm {POL}}}

\DeclareMathOperator{\Supp}{Supp}
\DeclareMathOperator*{\Res}{{Res}}

\newcommand {\e}{\varepsilon}

\newcommand {\set}[1]{\left\{ #1 \right\}}

\newcommand {\abs}[1]{\left|#1\right|}
\newcommand {\Norm}[1]{\left\|#1\right\|_2}
\newcommand {\absoo}[1]{\abs{#1}_\infty}
\newcommand {\K}{K}
\newcommand {\tr}{{\rm tr}}
\newcommand {\plv}{v}
\newcommand {\OOv}{\mathcal{O}_\plv}
\newcommand {\Kv}{\K_\plv}
\newcommand {\kv}{k} % KB: \bar{K} is bad: it is 1) not a "reduction of K",
\newcommand {\Gv}{G_\plv}

\begin{document}

\title{Factoring polynomials over global fields}

\author{Karim Belabas}
\address{Universit\'e Paris-Sud, D\'epartement de Math\'ematique, 91400
Orsay, France.}
\email{Karim.Belabas@math.u-psud.fr}

\author{Mark van Hoeij$\mbox{}^1$}\thanks{$\mbox{}^1$Supported by NSF grant 0098034}
\address{Florida State University, 211 Love building, Tallahassee, Florida
32306-3027, USA}
\email{hoeij@zeno.math.fsu.edu}

\author{J\"urgen Kl\"uners}
\address{Universit\"at Kassel, Fachbereich Mathematik und Informatik, 
Heinrich-Plett-Str. 40, 34132 Kassel, Germany.}
\email{klueners@mathematik.uni-kassel.de}

\author{Allan Steel}
\address{School of Mathematics and Statistics F07, 
University of Sydney NSW 2006, Australia}
\email{allan@maths.usyd.edu.au}

% \subjclass{Primary ; Secondary }

\begin{abstract} 
Let $\K$ be a global field and $f\in \K[X]$ be a polynomial.
We present an efficient algorithm which factors $f$ in
polynomial time.
\end{abstract}

\maketitle
\tableofcontents
 
%%%%%%%%%%%%%%%%%%%%%%%%%%%%%%%%%%%%%%%%%%%%%%%%%%%%%%%%%%%%%%%%%%%%%%%%%%%%

\section{Introduction}

  Let $\K$ be a global field. The goal of this paper is to present a practical
algorithm which factors polynomials $f\in \K[X]$ in polynomial time,
in particular for the cases $\K = \Q$ and $\K = \F_q(t)$. The seminal
Zassenhaus~\cite{Zassenhaus69} method to factor in $K[X]$ is as follows: we may
assume that $f$ is separable, integral and monic. First, compute a bound for the
factors of $f$, then find a non-archimedean place $\plv$ of $\K$ such that the
reduction $\bar{f}$ of $f$ modulo $\plv$ remains separable in $\kv[X]$
where $\kv$ is the residue field of $\plv$. Since $\kv$ is finite we
can factor $\bar{f}$ in $\kv[X]$ using well known algorithms.

Let $\Kv$ the completion of $\K$ at $\plv$. Let $\OOv$, resp.~$\OO$ be the
maximal order of $\Kv$, resp.~$\K$.
If $K = \Q$ then $\plv$ is a prime number, $\Kv$ the
$\plv$-adic numbers,
$\OOv$ the $\plv$-adic integers, $\OO = \Z$ and $k = \Z/v\Z$.

If $K = \F_q(t)$ then $\OO=\F_q[t]$ and we will choose a finite place $\plv$,
which corresponds to choosing an irreducible polynomial $\plv \in \F_q[t]$. If
$\alpha$ is a root of this polynomial, then $k \cong \F_q(\alpha)$,
$\OOv \cong k[[t-\alpha]]$ and
$\Kv \cong k((t-\alpha))$,
see also Section~\ref{sec:Fq[t]}.

After multiplying if necessary $f$ by an element of $\OO$, we may
assume $f \in \OO[X]$.
By Hensel's lemma, the factorization of $\bar{f}$ can be lifted to a
factorization
$$ f = \ell_f f_1 \cdots f_r $$
in $\Kv[X]$ where $\ell_f \in \K \subset \Kv$ is the leading coefficient
of $f$, and $f_1,\ldots,f_r$ are monic and irreducible in $\Kv[X]$.
We choose $\plv$ so that $f \in \OOv[X]$ and $\ell_f$ does not vanish mod $\plv$,
so $f_1,\ldots,f_r \in \OOv[X]$.
In actual computations, elements of $\OOv$ are computed modulo
$\plv^\ell$ for some $\ell > 0$ and lifted to $\OO$. For $a\in\OOv$
we write ``$a$ mod $\plv^\ell$'' for such a lift of $a$ to $\OO$.
This notation is extended to $\OOv[X]$ coefficientwise. By Hensel lifting the
irreducible factors of $\bar{f}$ we can
compute $f_1,\dots,f_r \mod \plv^\ell$ for any fixed
$\ell > 0$.

Let $g \in K[X] $ be a monic irreducible factor of $f$. Then 
$$g = f_1^{e_1} \cdots f_r^{e_r}$$
where $e_i\in\set{0,1}$ for all $1\leq i\leq r$.
If $\ell$ is large enough
compared to a bound on the coefficients of $g$, we may test for
given $e_1,\ldots,e_r \in \{0,1\}$ whether $f_1^{e_1} \cdots f_r^{e_r} \in K[X]$
by computing
$\ell_f f_1^{e_1} \cdots f_r^{e_r}$ mod $\plv^\ell$ and checking whether
this divides $f$ in $K[X]$.
This time, the lift ``$\ldots$ mod $\plv^\ell$'' to $\OO[X]$
is not arbitrary. Choosing the right lift is
straightforward if $\K$ is $\Q$ or $\F_q(t)$ since there are canonical minimal
lifts to $\OO$,
but requires care for general global fields (see \cite{Bel}
for the number field case).

The Zassenhaus algorithm finds the $e_i$ by an exhaustive enumeration, which
works very well if $r$ is small or the $K$-rational irreducible factors are
plentiful. Otherwise, we face combinatorial explosion and exponential behaviour.
The landmark paper by
Lenstra et al.~\cite{LLL} avoids this combinatorial problem by
constructing $K$-rational factors with lattice basis
reduction (LLL reduction).
The original paper assumes $K = \Q$, but was suitably generalized
by Arjen Lenstra~\cite{Len82} ($K$ a number field), then Pohst and
M\'endez~\cite{PoMe} ($K$ any global field). Unfortunately, although this
algorithm runs in polynomial time, it is rather slow in practice since its
worst case bounds require Hensel lift to huge accuracy, followed by the
LLL-reduction of correspondingly huge lattices. Mark van Hoeij~\cite{Hoe}
came back to the combinatorial problem and used a knapsack approach to solve
it for $K = \Q$, this was generalized to number fields by Belabas~\cite{Bel}.
These two papers stated no complexity bound. We shall describe a similar idea
over a general global field $K$, and show it runs in polynomial time,
although details will only be provided for the cases $K=\Q$ and $K=\F_q(t)$.
% It turns out to that this approach is highly practical. Although the algorithm
% is simpler if $K$ is a function field, we shall maintain a unified approach
% as long as possible.

\section{Notations}\label{notations}
Throughout the paper we will use the following notations:
Let $K$ be a global field of characteristic $p\geq 0$ with maximal order
$\OO$. We want to factor a separable polynomial $f\in \K[X]$ of
degree $n>1$. After multiplying by an element of $\K$ we may assume
that $f \in \OO[X]$.
Let $\plv$ be a non-archimedean place of $\OO$. We denote by
$\Kv$ the completion of $K$ at $\plv$, with maximal order 
$\OOv$, maximal ideal~$\plv$ and finite residue field $\kv$.
Let $\bar{f}$ be the image of $f$ in $\kv[X]$, and assume that $\bar{f}$
is still
separable. We also assume that the leading coefficient $\ell_f \in \OO$ of
$f$ does not vanish mod $\plv$, so the degree of $\bar{f}$ is still $n$.
% We remark that if $p>0$ then
% there may not be a place $\plv$ of degree 1 over $K$ that satisfies
% the above assumptions, and even if there is, it could still be better
% to use a place of degree $>1$. % (see Section~\ref{extending}).
In the number field case, instead of working with $\OO$ we can work
with a subring of $\OO$ if the computation of $\OO$ is too costly,
see \cite{Bel}.

We have the factorizations into irreducible elements
$$ f = \ell_f f_1 \cdots f_r\in\OOv[X],\quad
  \bar{f} = \bar{\ell}_f \bar{f}_1\cdots\bar{f}_r\in \kv[X], \quad\text{and}\quad
   f = \ell_f g_1 \cdots g_s\in \K[X].$$
% in $\OOv[X]$, $\kv[X]$ and $\K[X]$ respectively.
Furthermore,
$\ell_f g_i \in \OO[X]$.
Obviously $1 \leq s \leq r \leq n$.
We call the $f_i$ the \emph{local factors} and the $g_j$ the 
\emph{$K$-factors}.
We can not compute
$f_i \in \OOv[X]$ with infinite accuracy, but for any positive
integer $\ell$ we can
compute $f_i$ mod $\plv^\ell$, which is in $\OO[X]$.

\section{General description}
Our method relies on two main ideas:
\subsection{Linearize}
\label{3_1}

The logarithmic derivative is a group homomorphism from
the multiplicative group $\Kv(X)^*$ to the additive group $\Kv(X)$,
and has kernel $\Kv(X^p)^*$. The first main idea is to multiply this by $f$.
Then we obtain the following group homomorphism:
\begin{eqnarray*}
 \Phi: \Kv(X)^* / \Kv(X^p)^* & \rightarrow & \Kv(X) \\
           g & \mapsto & f g'/g.
\end{eqnarray*}
If $g$ is in the subgroup of $\Kv(X)^*/\Kv(X^p)^*$
generated by the local factors $f_i$, then $\Phi(g) \in \OOv[X]$.
If $g$ is in the subgroup generated by the $K$-factors $g_j$, then
$\Phi(g) \in \OO[X]$. To see this, take one such $g_j$.
Take any prime ideal of $\OO$ and let $w$ be the corresponding valuation
on $\K$, which is extended to a valuation on $\K[X]$ by taking $w(\sum c_i X^i) =
{\rm min}_i \ w(c_i)$.
Now $\Phi(g_j)$ is the product of $g_j'$ and $f/g_j$, both
of which are in $K[X]$,
but since $w(g_j') \geq w(g_j)$ we get
$w(f g_j'/g_j) \geq w(f g_j/g_j)$ which is $\geq 0$ since $f \in \OO[X]$.
So the valuation of $\Phi(g_j) \in K[X]$ is non-negative for any prime ideal
of $\OO$ and hence $\Phi(g_j) \in \OO[X]$.

Compared to the original algorithm of van Hoeij~\cite{Hoe}, we have replaced 
power sums by $f$ times the logarithmic derivative. To show the connection
we will define power sums.
Let $g \in K[X]$ be a monic
separable polynomial. Let $\alpha_1,\ldots,\alpha_m$ the zeros of $g$ in
an algebraic closure of $K$. For $j \geq 0$, the $j$'th power sum
($j$-th ``trace'') of $g$ is:
$$ {\rm Tr}_j(g) := \sum_{i=1}^m \alpha_i^j.$$
It is known that
$$ g'/g = \sum_{j\geq 0} {\rm Tr}_j(g) X^{-j-1},$$
which shows the relation between $g'/g$ and power sums.
Despite this relation,
our ``{\em $f$ times $g'/g$ approach}''  turns out to be more convenient
for complexity proofs than power sums,
and can also have practical advantages, particularly
when $f$ is not monic.

\subsection{Approximately solve knapsack}

Let $\Gv\subset\Kv(X)^*/\Kv(X^p)^*$ be the subgroup generated by the
local factors.
% and let $L := \Phi(\Gv)$ which is a $(\Z/p\Z)$-module.
Our goal is to find the subgroup $G\subset \Gv$
generated by the irreducible $K$-factors of $f$. To do this we first
construct the ``knapsack lattice'' $L$ in a similar way
as in \cite{Hoe}, except that instead of traces (power sums) of $f_j$
we use the coefficients of $\Phi(f_j)$.
We then reduce this lattice, which
means $\F_p$-linear Gaussian
elimination if $p > 0$, and LLL otherwise, for details see
Sections~\ref{4} and~\ref{sec:Fq[t]}.
Large basis vectors are then
discarded, yielding a sublattice $L'$ of $L$,
associated to a subgroup $G'$ of $\Gv$.

\subsection{Conclude}
\begin{theorem}\label{thm:bound}
We have $G = G'$ provided $\ell$ is large enough.
\end{theorem}
\begin{proof}
We will only sketch the proof, leaving the details to
Sections~\ref{4} and~\ref{sec:Fq[t]}.
% because we have not yet
% given details about the knapsack lattice $L$.
% We say that a factor appears in a rational function $g$
% if it divides either the numerator or
% the denominator of $g$.
If $G'$ is strictly larger than $G$, then by Lemma~\ref{specialg} below,
it contains an element, represented
by a rational function $g \in \Kv(X)^*$, $g \not\in \Kv(X^p)^*$,
% which we do not need to find explicitly,
such that
\begin{enumerate}
\item At least one $f_i$ divides $\Phi(g)$,
% not all local factors appear in $g$,
\item % None of the $g_j$ divide $\Phi(g)$, in fact,
None of the $\overline{g}_j$
divide $\overline{\Phi(g)}$ where the bar indicates reduction to $k[X]$.
% every $K$-factor has a local factor that appears in $g$,
\item $H := \Phi(g) \mod \plv^\ell$ is ``small''. \\
Indeed, it is a small perturbation
(see Lemma~\ref{specialg} below)
of a vector in the LLL basis of $L'$
that is small otherwise it would have been discarded.
\end{enumerate}
The clumsy argument for the third condition is
only needed if $p = 0$. If $p > 0$, all elements of $L'$ are small. Now, let
$g$ be as above,  $H := \Phi(g) \mod \plv^\ell$, and $R := \Res(f, H)$. Then
\begin{itemize}
\item $\Res(f, \Phi(g)) = 0$, hence $\plv^\ell \mid R$.  In fact,
$\plv^{\ell \sigma} \mid R$ where $\sigma$ is the sum of the degrees of
the $f_i$ that divide $\Phi(g)$. Item (1) above implies $\sigma>0$.
% local factor dividing $\Phi(g)$,
\item $R \neq 0$, because if $R$ was zero then $H$ would
be divisible by some $g_j$ so $\overline{H}$ would be divisible
by some $\overline{g}_j$, contradicting item (2) above.
% \item $\Res(f, H)$ is the determinant of the Sylvester matrix of $f$ and $H$
% and thus polynomially bounded in terms of the sizes of $f$ and $H$.
\end{itemize}
One obtains a contradiction if $\plv^\ell$ is larger than $R$. Since $R$
is  the determinant of the Sylvester matrix of $f$ and $H$, one obtains
a bound that is polynomial in terms of the sizes of $f$ and $H$.
\end{proof}

\begin{lemma}\label{specialg}
  Suppose $G \subsetneq G'$.
  Then there exists an element 
  $g \in G'\setminus G$ such that
  \begin{enumerate}
  \item $f_i \mid \Phi(g) \in\OOv[X]$ for some $1\leq i\leq r$.\label{en:1}
  \item $\overline{g}_j \nmid \overline{\Phi(g)}$ for all $1\leq j \leq s$.\label{en:2}
  \end{enumerate}
\end{lemma}
\begin{proof} Elements $g \in \Gv$ can be written in the form
$g = f_1^{e_1} \cdots f_r^{e_r} \cdot \Kv(X^p)^*$ where the
integers $e_i$ are defined mod $p$.
We view $e_i$ as element of $\Z/p\Z$, 
and then define the support of $g$ as
$\Supp g = \set{i \ | \ e_i\neq 0}$.
Since the $f_i$ are pairwise coprime and irreducible in $\OOv[X]$, we have
$$f_i \mid \Phi(f_j) \ \Longleftrightarrow \ i \neq j.$$
So $f_i \mid \Phi(g)$ iff $e_i$ is zero in $\Z/p\Z$, and $g_j \mid \Phi(g)$ iff
$\Supp g \bigcap \Supp g_j = \emptyset$.

The supports of $g_1,\ldots,g_s$
form a partition of $\set{1,\dots,r}$.
Choose any element $g \in G'\setminus G$. For all $1\leq j\leq s$ with
$\Supp g_j \cap \Supp g = \emptyset$, replace $g$ by $g_j g$.
Then condition \eqref{en:2} is satisfied (recall that $\overline{f}$ is separable),
and $g$ is still
in $G' \setminus G$. Write this $g$ as $f_1^{e_1} \cdots f_r^{e_r}
\cdot \Kv(X^p)^*$ with $e_i \in \Z/p\Z$.
Since $g$ is not in the group $G$ generated by $g_1,\ldots,g_s$,
there must be some $g_j$ for which $S_j := \{ e_i \ | \ i \in \Supp g_j \}$
contains more than one element. Then take an element $e \in S_j$
and replace $g$ by $g/g_j^e$.
Now $g$ satisfies both conditions \eqref{en:1} and \eqref{en:2}.
\end{proof} 
\begin{remark}
Given any $g \in G'\setminus G$, the above proof shows that
a ``small change'' suffices to obtain an element of $G'\setminus G$
that satisfies conditions (1) and (2).
\end{remark}

% If $p=0$ and if we take a ``short'' $g \in G'\setminus G$,
% then the construction in the above lemma results in a perturbation of $g$ that
% satisfies the two conditions, and that can be at most a small factor longer
% than the original $g$.

We have sketched a general proof and omitted the details.
Filling in these details is easy for the case $K = \F_q(t)$
discussed in Section~\ref{sec:Fq[t]}.
The details for $K = \Q$ require more work, which is
what we will do now.

\section{The case $K = \Q$} \label{4}
For $f \in \C[X]$ with leading coefficient $\ell_f$, let
$$ M(f) := \abs{\ell_f} \prod_{\abs{\alpha} > 1}
% _{\substack{f(\alpha) = 0\\ \abs{\alpha} > 1}}
\abs{\alpha}^{m_{\alpha}}$$
be the Mahler measure of $f$, where the product is taken over all
roots $\alpha \in \C$ of $f$ with absolute value $>1$, and $m_{\alpha}$
is the multiplicity of the root $\alpha$.
\begin{lemma} \label{logbound0}
If $f,g \in \C[X]$ and $g \mid f$ then
$$ \Phi(g) = \sum_{i=0}^{n-1} a_iX^i \in \C[X], \text{ with }
  \abs{a_i} \leq B_i := \binom{n-1}{i} n M(f). 
$$
\end{lemma}
\begin{proof}
We may assume $g$ is not a constant. Then the degree of $\Phi(g) \in \C[X]$ is $n-1$.
The Mahler measure of $\Phi(g)$ is bounded by $\deg(g) M(f)$ since
$M(A')\leq \deg(A) M(A)$, see \cite{Mah},
and $M(AB) = M(A)M(B)$ for any $A,B\in\C[X]$ \cite[p. 79]{MS}.
Bounding $\deg(g)$ by $n$,
the upper bound now follows by \cite[Lemma 2.1.9]{MS}.
\end{proof}

We restrict to the case $\K= \Q$, for a general number field follow \cite{Bel}.
We use the notation $\Z[X]_{<n}$ for all polynomials in $\Z[X]$ of degree $<n$.
We % identify $\Z[X]_{<n}$ with $\Z^n$ and
use $\Norm{.}$ for the $L^2$ norm on both $\Z^n$ and $\Z[X]_{<n}$.
\begin{corollary}\label{co:B}
With $f \in \Z[X]$ and $g$ any factor of $f$ in $\Q[X]$, we have
$\Phi(g) \in \Z[X]_{<n}$ and 
$$ \Norm{\Phi(g)} \leq B(f) := 2^{n-1}n \Norm{f}. $$
\end{corollary}
\begin{proof}
That $\Phi(g)$ is in $\OO[X]$ was proven in Section~\ref{3_1}.
Using Lemmata 2.1.8 and 2.1.9 in \cite{MS} we get that
 $\Norm{\Phi(g)} \leq 2^{n-1}M(\Phi(g))$.
 As in the proof of Lemma \ref{logbound0} we get that
  $2^{n-1}M(\Phi(g))\leq 2^{n-1}n M(f)$. Corollary
  2.1.5 in \cite{MS} states that $M(h)\leq \Norm{h}$ for all
  non constant  polynomials $h$ which finishes the proof.
\end{proof}

For $1 \leq j \leq s$ write the monic irreducible $K$-factors
as $g_j = f_1^{w_{j,1}} \cdots f_r^{w_{j,r}}$
with $w_{j,1},\ldots,w_{j,r} \in \{0,1\}$ and write $w_j := (w_{j,1},\ldots,w_{j,r})^\tr
\in \Z^r$ where $\tr$ denotes the transpose. Denote $W = \Z w_1 + \cdots + \Z w_s$.

If we have any basis $u_1,\ldots,u_s$ of $W$
then we can find $\{w_1,\ldots,w_s\}$ by computing the reduced echelon form
of $u_1,\ldots,u_s$, or by using the following shortcut: write $\{1,\ldots,r\}$
as a disjoint union of subsets in such a way that $i,j$
are in the same subset iff the $i$'th and $j$'th entry of $u$ are the same
for every $u$ in $u_1,\ldots,u_s$.

In the following let $I_r$ be the identity matrix of dimension $r$ and define
for $1\leq j \leq r$ the $a_{i,j}$ via
$$\Phi(f_j) {\rm \ mod \ } \plv^\ell = \sum_{i=0}^{n-1} a_{i,j} x^i.$$
Define the {\em all-coefficients} lattice $L$ as the span of the columns of the
following matrix: 
$$
A:=\left(\begin{matrix}I_r\\ A_1\\\end{matrix}\right), \mbox{ where }A_1:= 
\left(\begin{matrix}
a_{0,1} & \cdots & a_{0,r} \\
 \vdots & \ddots & \vdots      \\
a_{n-1,1}  & \cdots & a_{n-1,r} \\
\end{matrix}
\right).
$$
For $e=(e_1,\ldots,e_{r+n})^\tr \in L$ we denote the corresponding element 
$\Phi(f_1^{e_1}\cdots f_r^{e_r})$ of $\Phi(\Gv)$ as $\POL(e)$.
Each $K$-factor $g_j$ corresponds to a vector in $L$ we denote
by $\tilde{w}_j$, whose entries come from $w_j$ and $\Phi(g_j)$.
Then $\Norm{ \tilde{w}_j } \leq \sqrt{
\Norm{w_j}^2 + B^2} \leq B' := \sqrt{r^2+B^2}$ where $B = B(f)$
is as in Corollary~\ref{co:B}.
\begin{theorem} \label{ThmG}
Let $f \in \Z[X]$ separable and $B'$ as above.
% Let $\ell$ be any integer such that 
% \begin{equation}\label{eq:BQ}
% \plv^{\ell/n} > \Norm{f} (2^{n-1} + n) B' (1+B').
% \end{equation}
Let $b_1,\dots,b_n$ an LLL-reduced basis for the all-coefficients lattice $L$,
let $b^*_1,\dots,b^*_n$ the associated Gram-Schmidt
orthogonalized basis, and let $t$ the smallest index such that $\Norm{b^*_j}
> B'$ for all $j > t$.
% Then $t = s$ and the projection of
% $b_1,\dots,b_s$ on the first $r$ entries is a $\Z$-basis of $W$.
Let $L' := \Z b_1 + \cdots + \Z b_t$.
If
\begin{equation}\label{eq:BQ}
\plv^{\ell/n} > \Norm{f} (2^{n-1} + n) B' (1+B').
\end{equation}
then the projection of $L'$ on the first $r$ entries is $W$.
\end{theorem}
\begin{proof}
It follows from the proof of (1.11) in \cite{LLL} that every $w \in L$
with $\Norm{w} \leq B'$ is in $L'$.
So $\tilde{w}_1,\ldots,\tilde{w}_s \in L'$
and hence the projection of $L'$ on the first $r$ entries
contains $W$. Assume it is strictly larger,
then $\POL(b_u) \not\in \Phi(G)$ for some $1\leq u\leq t$. From the
properties of LLL-reduced bases, $\Norm{b_u} \leq \gamma^{t - 1} B'$, where
$\gamma > 4/3$ is a number that can be chosen in the reduction algorithm
(one may set $\gamma := 2$ as in the original LLL paper). More precisely,
$\Norm{b_u^*} \leq \gamma^{t - u} \Norm{b^*_t}$ and $\Norm{b_u}\leq
\gamma^{u-1}\Norm{b^*_u}$ for all $u \leq t \leq n$.

Using Lemma~\ref{specialg} as in the proof of
Theorem~\ref{thm:bound} one can show that
there exists a vector $g \in L'$ such that $f_i \mid \POL(g)$
for some $1\leq i\leq r$, and $\plv^\ell \mid
\Res(f, H) \neq 0$, where $H := \POL(g) \mod \plv^\ell$. From the proof of
Lemma~\ref{specialg}, this vector
$g$ may be obtained by first adding a subset of $\{\tilde{w}_1,
\ldots,\tilde{w}_s\}$
to $b_u$, yielding a vector $b$ such that
$$\Norm{b} \leq (\gamma^{t-1} + s) B',$$
and then by adding to $b$ a vector of the form $e \tilde{w}_i$
for some integer $e$
with $\abs{e} \leq
\absoo{b} \leq \Norm{b}$. Hence
$$\Norm{H} \leq \Norm{g} \leq (\gamma^{t-1} + s) B'(1 + B').$$
 From the preceding discussion and Hadamard's bound,
$$ \plv^\ell \leq \abs{\Res(f, H)} \leq \Norm{f}^{\deg{H}} \Norm{H}^n, $$
and we may bound $\deg H, s, t \leq n$ to derive a contradiction with
\eqref{eq:BQ}.
\end{proof}

 From this theorem, we obtain $\ell\log \plv = O(n^2 + n\log \Norm{f})$.
Since Hensel lifting and lattice reduction are polynomial time algorithms,
we see that $W$ can be computed in polynomial time.

Although there is no practical reason for doing so (since power sums
do not offer advantages over coefficients of $\Phi$), one could now
use the relation between power sums (called traces in \cite{Hoe})
and $\Phi$ to show the algorithm in \cite{Hoe} is polynomial
time provided that one uses what we call the {\em all-traces} version
of the algorithm. This version uses all of the traces numbered $1,\ldots,n-1$
at the same time, so the lattice reduction takes place in $\Z^{r+n-1}$.
% (the $0$'th trace and the $n-1$'th coefficient of $\Phi$ can be ignored
% in the implementation).
 From a practical point of view, the all-traces
and all-coefficients versions are slow and thus not interesting.

The main question is whether practical versions of the algorithm
run in polynomial time.
% This question is only well defined
% if one specifies which recipe the algorithm must use for
% choosing the lattice.
Using one trace at a time works very well in practice, see \cite{Bel}.
We will show that the ``one coefficient at a time'' version
factors in $\Q[X]$ in polynomial time (the same must then also be
true for one trace at a time).
% This version should have the
% same running times in practice as those reported in \cite{Bel}
% if we start with the $X^{n-2}$ coefficient.
  % The relation between traces and $\Phi$ then implies that the algorithm
  % itself (with traces) also factors in $\Q[X]$ in polynomial time.
  % Note that for implementations we propose to use coefficients of $\Phi(f_i)$
  % instead of traces, although we do not expect much difference in
  % running times.

Let $B_i$ be the bound for $\abs{a_i}$ given
in Lemma~\ref{logbound0}.
For $0 \leq i \leq n-1$ and $g \in \Kv(X)^*$ write $T'_i(g) \in \Z$
the coefficient of $X^i$ in $\Phi(g)$ mod $\plv^\ell$.
Let $T_i(g) := T'_i(g)/B_i \in \Q$.
Now Lemma~\ref{logbound0} says that if $g$ is a $K$-factor of $f$, then
$\abs{T_i(g)} \leq 1$.

\begin{proposition}
\label{sequenceL}
One can compute a sequence of lattices $L_{n-1},L_{n-2},\ldots,L_{0}$ with
the following properties:
\begin{enumerate}
\item $\Z^r = L_{n-1} \supseteq L_{n-2} \cdots \supseteq L_{0} \supseteq W$
\item $L_i = \Z b_{i,1} + \cdots + \Z b_{i,r_i}$ for some integer $r_i$
and some vectors $b_{i,j} \in \Z^r$ with the following properties:
\begin{enumerate}
\item $\Norm{ b_{i,j} } \leq (r+2) \gamma^r$.
\item If $b_{i,j} = (e_1,\ldots,e_r)^\tr$ then
$T_i(f_1^{e_1} \cdots f_r^{e_r}) \leq (r+2) \gamma^r$
\end{enumerate}
\end{enumerate}
where % $\gamma$ can be taken as $\gamma=2$.
$\gamma > 4/3$ is a number that can be chosen in the reduction algorithm
(one may set $\gamma := 2$ as in the original LLL paper).
\end{proposition}
\begin{proof}
If $i=n-1$ we may take $b_{i,1},\ldots,b_{i,r_i}$ as the standard basis of $\Z^r$.
If $i<n-1$ then we may assume that $L_{i+1} = \Z b_{i+1,1} + \cdots + \Z
b_{i+1,r_{i+1}}$ has been computed and define $b'_j$ as follows:
First write $b_{i+1,j} = (e_1,\ldots,e_r)^\tr$, then compute
$a := e_1 T_i(f_1) + \cdots + e_r T_i(f_r)$ and set
$b'_j := (e_1,\ldots,e_r,a)^\tr \in \Z^r \times \Q$.
Now let $L' := \Z b'_1 + \cdots + \Z b'_{r_{i+1}} + \Z P$
where $P = (0,\ldots,0,\plv^\ell/B_i)^\tr$.
Let $b_1,b_2,\ldots$ be an LLL-reduced basis of $L'$, let
$b_1^*, b_2^*, \ldots$ the associated orthogonalized basis, and let $r_i$
be the smallest index such that
$\Norm{b^*_j} > r+2$ for all $j > r_i$.
Now define $b_{i,j}$ as the projection of $b_j$ on the first $r$ entries
and let $L_i := \Z b_{i,1} + \cdots + \Z b_{i,r_i}$.

Consider the vector $w_j$ corresponding to the $K$-factor $g_j$
and let $w'_j$ be the corresponding vector in $L'$.
The first $r$ entries of $w'_j$
are in $\{0,1\}$, and the last entry equals $T_i(g_j) \in \Q$
which has absolute value $\leq 1$ by
Lemma~\ref{logbound0}. Hence, $\Norm{w'_j} \leq \sqrt{r+1} < r+2$.
Then it follows from the proof of (1.11) in \cite{LLL} that $w_j \in L_i$
and hence $W \subseteq L_i$.
By the properties of an LLL-reduced basis, we have
$\Norm{ b_j } \leq (r+2) \gamma^r$ when $j \leq r_i$ which implies
(2a) resp. (2b) since projecting on the first $r$ entries
resp. last entry does not make a vector longer.

The lattice $L'$ to be reduced was in
$\Z^r \times \Q$. Lattice reduction in $\Z^{r+1}$ is
more efficient, so we round each of the
numbers $T_i(f_1),\ldots,T_i(f_r),\plv^\ell/B_i$ to the nearest integer.
Then we obtain a lattice $L' \subseteq \Z^{r+1}$ but now we have
introduced rounding errors.
Consider again the vectors $w_j \in W$ and $w'_j \in L'$.
If $w_j$ has $\sigma$ entries equal to 1, then the last entry of $w'_j$
is the sum of $\sigma$ of elements of
$\{T_i(f_1),\ldots,T_i(f_r)\}$ plus an integer in the interval
$(-\sigma/2,\sigma/2)$
times $\plv^\ell/B_i$. We introduced an error $\leq 0.5$ in each of the
numbers $T_i(f_1),\ldots,T_i(f_r),\plv^\ell/B_i$.
Then the total rounding error in the last entry of $w'_j$ is
less than $0.5 (\sigma + \sigma/2)$ which is less than $r$,
so this entry will have absolute value $< r+1$.
Then $\Norm{w'_j} < \sqrt{\sigma + (r+1)^2} < r+2$.
The proposition is stated with $r+2$ instead
of $\sqrt{r+1}$ so that the bounds can still be used for
practical implementations that round
$T_i(f_1),\ldots,T_i(f_r),\plv^\ell/B_i$ to $\Z$.
\end{proof}

If $L_{i+1}$ is known, then the computation
of $L_i$ in the proposition involves a lattice reduction
in $\Z^{r+1}$ of a lattice with determinant $\plv^\ell / B_i$
(rounded to the nearest integer).
If $\plv^\ell / B_i$ is large, then we get
a big practical improvement by doing this lattice reduction incrementally
in the way it is described in Section 2.4 in
\cite{Bel}, reducing one large-determinant
lattice reduction to a sequence of smaller lattice reductions that at
the end produce the same result.

\begin{lemma}
\label{lemr}
With the notations of Proposition~\ref{sequenceL},
the following holds for every $n-1 \geq i \geq i' \geq 0$.
If $e=(e_1,\ldots,e_r)^\tr$ is an element of $\{b_{i',1},\ldots,b_{i',r_{i'}}\}$
then
\[ T_{i}(f_1^{e_1} \cdots f_r^{e_r}) \leq 2^{O(r^2)} \]
\end{lemma}
\begin{proof}
The entries of the $b_{i,j}$ and $e$
are bounded by $(r+2)\gamma^r = 2^{O(r)}$.
Since $e \in L_{i'} \subseteq L_i$ we
can write $e = \sum_{j=1}^{r_i} c_j b_{i,j}$ for some $c_j \in \Z$
that can be found by solving linear equations. With Cramer's rule
one finds $\abs{c_j} \leq 2^{O(r^2)}$. Multiplying this by $r_i$
and by the bound given in (2b) in Proposition~\ref{sequenceL} one
obtains the bound $2^{O(r^2)}$.
\end{proof}

% To prove that the ``one trace at a time'' approach (or actually:
% To prove that the ``one coefficient of $\Phi$ at a time'' version
% is polynomial time, we must show that $L_0 = W$ for some polynomially
% bounded $\ell$. First a lemma.

\begin{theorem}
$L_0 = W$ for some $\ell$ with $\ell \log \plv$ polynomially
bounded in terms of the degree of $f$ and $\log \Norm{f}$.
\end{theorem}
\begin{proof}
If $L_0 \neq W$ then let $e$ be one of the vectors
$b_{0,j}$ from Proposition~\ref{sequenceL} that is not in $W$.
Write $e=(e_1,\ldots,e_r)^\tr$ and $g = f_1^{e_1} \cdots f_r^{e_r}$.
Write $\Phi(g) = \sum c_i X^i$. Then the corresponding
vector in the all-coefficients lattice (see Theorem~\ref{ThmG})
is $\tilde{e} := (e_1,\ldots,e_r,c_0,\ldots,c_{n-1})^\tr$
where $c_0,\ldots,c_{n-1}$ are bounded in absolute value by $2^{O(r^2)}$
by Lemma~\ref{lemr}.
Applying the process in the proof of Lemma~\ref{specialg}
we obtain a new vector $e'$ whose length differs
at most $(s + {\rm max}\{e_1,\ldots,e_r\})B'$ from $e$.
The last $n$ entries of this vector are the coefficients of
a polynomial $H \in \Z[X]_{<n}$ and
% These coefficients are polynomially
% bounded, and
we have $\plv^\ell \mid {\rm Res}(f,H) \neq 0$ in the same
way as in Theorem~\ref{ThmG}. This implies that $\log \plv^\ell$
is polynomially bounded.
\end{proof}

We propose to implement the ``one coefficient at a time'' approach
in the following way: start with a value for $\ell$ that
is {\em at most} as large as what one would use in the Zassenhaus
approach.
Then, compute $L_{n-1}, L_{n-2}, \ldots$ until we find $W$.
If we reach $L_0$ and we still have not found $W$ then we must
increase $\ell$.
The computation of each $L_i$ should be done using the incremental
strategy of Section 2.4 in \cite{Bel}.
Then one has a polynomial time algorithm that runs very well in
practice, with running times that are essentially the same
as those reported in \cite{Bel} for $\Q[X]$.

\section{The case $K=\F_q(t)$} \label{sec:Fq[t]}

Now $\OO=\F_q[t]$,
and the place $\plv$ corresponds to an irreducible polynomial
in $\F_q[t]$, which we shall also denote as $\plv$.
Let $f \in \OO[X]$. We want to factor $f$, viewed as element of $\F_q(t)[X]$. We assume
that $f$ is separable.
Denote $\alpha$ as a root of $\plv \in \F_q[t]$, then the residue field $\kv = \F_q[t]/(\plv)$ is
isomorphic to $\F_q(\alpha)$.
We choose $\plv$ in such a way that $\bar{f}$, the image of $f$ in $\kv[X]$, is squarefree
and of the same $X$-degree as $f$.
We get the factorization
$$\bar{f} = \bar{\ell_f} \bar{f}_1\cdots \bar{f}_r \in \kv[X].$$
Representing
$t-\alpha$ with a new variable
$\tilde{t}$, the map $t \mapsto \tilde{t}+\alpha$ is an isomorphism from
$\F_q[t]/(\plv^\ell)$ to $\F_q(\alpha)[\tilde{t}]/(\tilde{t}^\ell)$.
Taking limits, one finds an isomorphism from
\[ \OOv = \lim_{\leftarrow} \F_q[t]/(\plv^\ell)
\]
to
\[ \F_q(\alpha)[[\tilde{t}]] = \lim_{\leftarrow} \F_q(\alpha)[\tilde{t}]/(\tilde{t}^{\ell}).
\]
By Hensel's lemma, we get a factorization
$$f = \ell_f f_1 \cdots f_r \in \OOv[X].$$
If $g \in \OOv[X]$ we denote ``$g$ mod $\plv^\ell$'' as the unique lift
of $g$ to $\F_q[t,X]$ whose $t$-degree is smaller than the $t$-degree
of $\plv^\ell$.
We can not compute $f_i \in \OOv$ with infinite
accuracy, however, for any integer $\ell>0$
we can compute $f_i$ mod $\plv^\ell$,
which is an element of $\F_q[t,X]$.

Note that the above technicalities with $\OOv$ become easier if we take
$\plv = t$ so that $\tilde{t}=t$.
However, we can not always do this, we can
only take $\plv=t$ if $f(t=0,X)$ is
square-free and of the same degree as $f$.
% First we need a bound.
\begin{lemma}\label{logboundp}
  Let $g\in \F_q[t][X]$ be a polynomial which divides $f$ then 
$$\Phi(g) = \sum_{i=0}^{n-1} a_i(t) x^i \in \F_q[t][x] \mbox{ with }
  \deg(a_i) \leq B_i := \deg_t(f),$$
  where $\deg_t$ denotes the $t$-degree.
\end{lemma}
\begin{proof}
  From $\Phi(g)=f g' / g$
  we get $\deg_t(\Phi(g)) + \deg_t(g) =
  \deg_t(g') + \deg_t(f)$.
  Since $\deg_t( g' )
  \leq \deg_t(g)$ we get the wanted bound.
\end{proof}

The idea is as follows. Let $g \in \Gv$.
If the degree of one of the coefficients of $\Phi(g)$ mod $\plv^\ell$
exceeds the degree bound $B_i$
then $g$ is not a $K$-factor
of $f$. We use this to replace the Zassenhaus combinatorial
search by linear algebra.

As in the rational case we introduce the lattice $W \subseteq (\Z/p\Z)^r$ generated by the 
exponent vectors of the monic irreducible factors $g_1,\ldots,g_s$ of $f$ in $\F_q(t)[X]$.
Let $L$ be some subspace of $(\Z/p\Z)^r$ that contains $W$.
We start with $L = (\Z/p\Z)^r$.
For $e=(e_1,\ldots,e_r)^\tr \in L$ we denote by $\POL(e)$ the polynomial
$$\Phi(f_1^{e_1} \cdots f_r^{e_r}) {\rm \ mod \ } \plv^{\ell}.$$

Our goal is to compute a subspace of $L' \subseteq L$ which still contains $W$.
Write
$$\POL(\e_j)=\sum_{i=0}^{n-1} a_{i,j} X^i  \ \ \ \ (1\leq j\leq r)$$
where $\e_1,\ldots,\e_r$ is the standard basis of $(\Z/p\Z)^r$.
% In case that $\e_j$ corresponds to a true factor, Lemma~\ref{logboundp}
% shows that $\deg(a_{i,j})\leq B_i$.
Let $m_i = B_i+1$ and let $\sigma$ be the $t$-degree of $\plv^\ell$.
We define
% $$\phi_{m_i}: \F_q[t] \rightarrow \F_q^{\sigma-m_i} \ \ \ \ 
$$ \phi_{m_i}(\sum_k c_k t^k)
:= (c_{m_i},\ldots,c_{\sigma-1})^\tr$$ 
and
\[
A_i := 
\left( \phi_{m_i}(a_{i,1}) \cdots \phi_{m_i}(a_{i,r}) \right)
% \left(
% \begin{matrix}
% b_{1,1} & \dots  & b_{1,\sigma-m_i} \\
% \vdots       & \ddots & \vdots       \\
% b_{r,1} & \dots  & b_{r,\sigma-m_i} \\
% \end{matrix}
% \right):=
% \left(
% \begin{matrix}
% \phi_{m_i}(a_{i,1})  \\
% \vdots     \\
% \phi_{m_i}(a_{i,r})  \\
% \end{matrix}
% \right) 
% \in \F_q^{\, \sigma-m_i \times r}
\]
which is an $(\sigma-m_i) \times r$ matrix with entries in $\F_q$,
and $A_i e=0$ for all $e \in W$.
% Therefore a kernel computation yields a subspace of $W'$
% containing $W$.
Now $L$ and $W$ are subspaces of $\F_p^{\,r}$ and $A_i$ is defined
over $\F_q$.
For $q=p^w$ write $\F_q=\F_p \gamma_1
 + \cdots \F_p \gamma_w$ and define
$$\psi: \F_q \rightarrow \F_p^{\,w}, \ \ \ \
\sum_{l=1}^{w} c_l \gamma_l \mapsto
(c_1,\ldots,c_{w})^\tr,$$
where $\tr$ denotes the transpose.
We define $\tilde A_i$ as follows: replace every entry $c$ of $A_i$
by $\psi(c)$. Since $\psi(c)$ is a column vector (because of the transpose in
its definition) with $w$ entries we see that $\tilde A_i$ is an $w(\sigma-m_i)
\times r$
matrix with entries in $\F_p$.
% \left(
% \begin{matrix}
% \psi(b_{1,1}) & \dots  & \psi(b_{1,\sigma-m_i}) \\
% \vdots       & \ddots & \vdots       \\
% \psi(b_{r,1}) & \dots  & \psi(b_{r,\sigma-m_i}) \\
% \end{matrix}
% \right)\in \F_p^{\, rw \times \sigma-m_i}
% \]
We still have $\tilde A_i e = 0$ for all $e \in W$.
% Therefore a kernel computation over $\F_p$
% yields a subspace $L'$ of $L$ still containing $W$. 
Now let $L'$ be the intersection of the kernels of
$\tilde{A}_0,\ldots,\tilde{A}_{n-1}$.
Then $L'$ contains $W$.

Let $B$ be a degree bound which can be easily computed using Theorem
\ref{thm:bound}.  E.g. we can take
$B=(2n-1)\deg_t(f)$ when we use the estimate
from Lemma \ref{logboundp} and the properties of the Sylvester matrix.
% We remark that these bounds can be easily improved.
% In case that $\ell>B+1$ and we apply this procedure to all
% $0\leq i \leq n-1$.
Theorem \ref{thm:bound} guarantees that $L'$ will be $W$ when $\sigma$
(the $t$-degree of $\plv^\ell$) is larger than $B$.
% We already get this result if we only use the $\F_q$-kernel of the
% matrices $A_i$.
Altogether we have proved
\begin{theorem}
  If the $t$-degree of $\plv^\ell$ is larger than $(2n-1)\deg_t(f)$ and $L'$
  is the intersection of the kernels of $\tilde{A}_i$, $i=0,\ldots,n-1$
  then $L' = W$. This leads to an algorithm that produces the
  factorization of a separable polynomial $f\in \F_q[t][X]$ in polynomial time.
\end{theorem}
\begin{remark}
\label{biv}
If the total degree of $f$ as bivariate polynomial is $n$,
then one can replace $B_i := {\rm deg}_t(f)$ in Lemma \ref{logboundp}
by $B_i := n-1-i$. The proof is essentially the same, except that
the degree w.r.t. $t$ should be replaced by the total degree.
Then we can replace $(2n-1)\deg_t(f)$ by $n(n-1)$ in the above theorem.
\end{remark}
\begin{remark}
In \cite{issac04} the authors followed our ``$f$ times $g'/g$''
approach found in a previous version of this paper and were able to
improve the quadratic bound $n(n-1)$ to a linear bound when $p>n(n-1)$.
\end{remark}
Note that in an implementation, one would start with a small value
for $\ell$, increasing $\ell$ as long as $L'$ is not $W$.
To improve practical performance, we can replace
the bounds $B_i$ from Lemma \ref{logboundp} or Remark~\ref{biv}
by the sharper bound given in the lemma below.

Denote $N(f) \subset \R^2$ as the Newton polygon of $f$, which
is defined as the convex hull of all points $(i,j)$ for which
the coefficient of $t^i X^j$ in $f$ is non-zero.
If $S_1,S_2 \subset \R^2$ 
then define $S_1+S_2 := \{s_1+s_2 \ | \ s_1 \in S_1,
s_2 \in S_2\}$.
\begin{lemma}
% Let $C_i$ be the largest real number for which $(C_i, i) \in
% N(f) +  \{ (0,-1) \}$ and let $B_i$ be $C_i$ rounded down to
% an integer.
Let $B_i := {\rm sup} \{ m \in \N \ | \ (m,i) \in N(f) + \{(0,-1)\} \ \}$.
Let $g\in \F_q[t][X]$ be a polynomial which divides $f$ then
$$\Phi(g) = \sum_{i=0}^{n-1} a_i(t) x^i \in \F_q[t][x] \mbox{ with }
  \deg(a_i) \leq B_i.$$
\end{lemma}
\begin{proof}
It is well known that $N(g h)  
= N(g) + N(h)$ for all $g,h \in \F_q[t,X]$.
It is also clear
that $N(g') \subseteq N(g) + \{ (0,-1) \}$.
Then $N( \Phi(g) ) = N(f/g \cdot g') = N(f/g) + N(g')
\subseteq N(f/g) +  N(g) + \{ (0,-1) \} = N(f) +  \{ (0,-1) \}$.
\end{proof}

%%%%%%%%%%%%%%%%%%%%%%%%%%%%%%%%%%%%%%%%%%%%%%%%%%%%%%%%%%%%%%%%%%%%%%%%%%
%%%%%%%%%%%%%%%%%%%%%%%%%%%%%%%%%%%%%%%%%%%%%%%%%%%%%%%%%%%%%%%%%%%%%%%%%%

\newcommand{\etalchar}[1]{$^{#1}$}

\end{document}